\documentclass[11pt]{article}

\begin{document}
 
\newcommand{\e}{\epsilon}
\renewcommand{\a}{\alpha}
\renewcommand{\b}{\beta}
\newcommand{\om}{\omega}
\newcommand{\La}{\Lambda}
\newcommand{\la}{\lambda}
\newcommand{\p}{\partial}
\newcommand{\de}{\delta}
\newcommand{\De}{\Delta}
\newcommand{\id}{\int_D}
\newcommand{\na}{\nabla} 
 
\newtheorem{corollary}{Corollary} 
\newtheorem{proposition}{Proposition} 
\newtheorem{lemma}{Lemma}  
\newtheorem{theorem}{Theorem} 
 
\def\sqg{G(\sqrt{(x-\xi)^2+(y-\eta)^2})} 
\def\kernel{G(\sqrt{x^2+y^2})} 
\def\be{\begin{eqnarray}} 
\def\ee{\end{eqnarray}} 
\def\qed{\hbox{\vrule width 6pt height 6pt depth 0pt}}

\title{Asymptotic Dynamical Difference between the Nonlocal and Local 
	Swift-Hohenberg Models 
\footnote{Author for correspondence: Professor Jinqiao Duan,
Department of Mathematical Sciences,  
Clemson University, Clemson, South Carolina 29634, USA.
$\;$ E-mail:  duan@math.clemson.edu; $\;$ Fax: (864)656-5230.} }

\author{Guoguang Lin$^1$, Hongjun Gao$^2$,  
Jinqiao Duan$^3$ and Vincent J. Ervin$^3$  
   \\ 
   \\ 
\\1. Graduate School, Chinese Academy of Engineering Physics \\  
P. O. BOX 2101, Beijing  100088, and
Department of Mathematics \\  
Yunnan University, Kunming  650091,  China. 
\\ 
\\2. Laboratory of Computational Physics\\ 
Institute of Applied Physics and Computational Mathematics\\ 
Beijing, 100088, China. 
\\ 
\\3. Department of Mathematical Sciences\\  
Clemson University, Clemson, South Carolina 29634, USA.  }

\date{June 11, 1998 } 
 
\maketitle 
 
\begin{abstract} 
 
In this paper the difference in the asymptotic dynamics between
the nonlocal and local two-dimensional Swift-Hohenberg models
is investigated. It is shown that the bounds  for the dimensions 
of the global attractors for the nonlocal 
and local Swift-Hohenberg models differ by an absolute constant, 
which depends only on the Rayleigh number, 
and upper and lower bounds of the kernel of the nonlocal
nonlinearity.  Even when this kernel of the nonlocal operator
is a constant function, the dimension bounds  of the global attractors
still differ by an absolute constant depending on the Rayleigh number.

\bigskip	  
\bigskip  
{\em Running Title:} Nonlocal Swift-Hohenberg Model
 
 
	  
\end{abstract}

\newpage

   \section{Introduction} 
 
Fluid convection due to density gradients arises in 
geophysical fluid flows in the atmosphere, oceans
and the earth's mantle. 
The Rayleigh-Benard convection is a prototypical model 
for fluid convection, aiming at predicting spatio-temporal 
convection patterns.  
The mathematical model for the 
Rayleigh-Benard convection involves nonlinear Navier-Stokes  
partial differential equations coupled with the temperature 
equation. When the Rayleigh number is near the onset of the 
convection, the Rayleigh-Benard convection model may be approximately 
reduced to an amplitude or order parameter equation, as derived 
by Swift and Hohenberg (\cite{Swift}).  

In the current literature, most work on the 
Swift-Hohenberg model deals with the following one-dimensional equation  
for $w(x,t)$, which is a  
localized, one-dimensionalized version of the model
originally derived by Swift and Hohenberg (\cite{Swift}),
\begin{eqnarray} 
	w_t & = & \mu w -(1+\p_{xx})^2 w - w^3. 
\label{sh2} 
\end{eqnarray}
The cubic term $w^3$ is used as an approximation of
a nonlocal integral term. 
For the (local) one-dimensional Swift-Hohenberg 
equation (\ref{sh2}), there has been some recent research on  
propagating or steady patterns (e.g., \cite{Eckmann_Wayne}, \cite{Hilali}, 
\cite{Lerman}).     
Mielke and Schneider(\cite{Mielke})  proved the existence of   
the global attractor 
in a weighted Sobolev space on the whole real line.  
Hsieh et al. (\cite{Hsieh}, \cite{Hsieh2}) 
remarked that the elemental instability mechanism is 
the negative diffusion term $-w_{xx}$.

Roberts (\cite{Roberts}, \cite{Roberts2}) 
recently re-examined  the rationale for using the  
Swift-Hohenberg model as a reliable model of the spatial 
pattern evolution in specific physical systems. 
He   argued  that,  
although the localization 
approximation used in (\ref{sh2}) 
makes some sense in the one-dimensional case,     
this approximation is deficient  
in the two-dimensional convection problem and 
one should use the nonlocal Swift-Hohenberg model 
(\cite{Swift}, \cite{Roberts}, \cite{Roberts2}): 
\begin{eqnarray} 
 u_t  = \mu u -(1+\De)^2 u  
	- u \int_D G(\sqrt{(x-\xi)^2+(y-\eta)^2}) u^2(\xi, \eta, t)d\xi d\eta , 
\label{sh} 
\end{eqnarray} 
where $u=u(x,y,t)$ is the unknown amplitude function,
$\mu$ measures the difference of the Rayleigh number
from its critical onset value,
$\De = \p_{xx} + \p_{yy}$ is the Laplace operator, and 
$G(r)$ is a given radially symmetric function ($r=\sqrt{x^2+y^2}$). 
The equation is defined for
$t>0$ and $(x,y) \in D$, where $D$ is a bounded planar domain with 
smooth boundary $\p D$. 
 
The two-dimensional version of the local  Swift-Hohenberg equation  
for $u(x, y, t)$ is 
\begin{eqnarray} 
	u_t & = & \mu u -(1+\De)^2 u - u^3. 
\label{sh3} 
\end{eqnarray} 
Here $u^3$ is used to approximate the nonlocal term in (\ref{sh}). 
 
Roberts (\cite{Roberts}, \cite{Roberts2}) noted that the range of Fourier
harmonics generated by the nonlinearities is fundamentally
different in two-dimensions than in one-dimension. This difference
requires a more sophisticated treatment of two-dimensional
convection problem, which leads to nonlocal nonlinearity
in the Swift-Hohenberg model. He also argued
that nonlocal operators naturally appear in systematic derivation of
simplified models for pattern evolution, and nonlocal operators
also permit symmetries which are consisitent with
physical considerations.

In this paper, we discuss the difference between
nonlocal and local two-dimensional 
Swift-Hohenberg models (\ref{sh}), (\ref{sh3}),
from a viewpoint of asymptotic dynamics. 
We show that the bounds  for the dimensions 
of the global  attractors for the nonlocal 
and  local Swift-Hohenberg models differ by an absolute 
constant, which depends only on
the the Rayleigh number, and upper and 
lower bounds of the kernel of the nonlocal
nonlinearity.  Even when this kernel is a constant function, 
the dimension bounds  of the global attractors
still differ by a constant depending on the Rayleigh number.
In \S 2 and \S 3, we will consider the 
nonlocal and local Swift-Hohenberg models, respectively.
Finally in \S 4, we summarize the results.

\section{Nonlocal Swift-Hohenberg Model}

In this section, we discuss the global attractor 
and its dimension estimate for the nonlocal
Swift-Hohenberg model (\ref{sh}).
In the following we use the abbreviations $L^2=L^2(D)$, 
$L^{\infty}=L^{\infty}(D)$,  $H^k = H^k(D)$ and $H^k_0 =H^k_0(D)$  ($k$ is 
a non-negative integer) 
for the standard Sobolev spaces.  
Let $(\cdot, \cdot)$, 
$\| \cdot \| \equiv \| \cdot \|_2$ denote the standard  
inner product and norm in $L^2$, 
respectively. The norm for $H^k_0$  is  
$\| \cdot \|_{H^k_0}$. Due to the Poincar\'e inequality,  
$\|D^k u\|$ is an equivalent norm in $ H^k_0$.

We rewrite the two-dimensional nonlocal Swift-Hohenberg equation (\ref{sh}) 
as 
\begin{eqnarray} 
u_t + \a u +2\De  u + \De^2 u 
  + u \int_D G(\sqrt{(x-\xi)^2+(y-\eta)^2}) u^2(\xi, \eta, t)d\xi d\eta=0, 
\label{sh4} 
\end{eqnarray} 
where  $\a = 1 - \mu$. 
This equation is supplemented with the initial condition 
 \be 
 u(x,y,0)=u_0(x,y), 
	\label{sh5} 
 \ee  
and the  boundary conditions 
 \be 
 u|_{\p D}   =    0, \;\;  \frac{\p u}{\p n} |_{\p D} = 0,   
	\label{sh6} 
 \ee 
where $n$ denotes the unit outward normal vector of the boundary $\p D$. 
 
In this paper, we assume the following conditions  
for every $t \geq 0 $ and $(x,y) \in D$,   
 \be
  0 < b  \leq  \kernel  \leq  a,  \; \mbox{and}  \;
 G, \; \na G, \; \De G \in L^{\infty}(D),
 \label{sh7}
 \ee
where $a, b > 0$ are some positive constants
and $\na = (\p_x, \p_y)$ is the gradient operator. 
Denote $K_1= \|\na G\|_{\infty}$ and $K_2=\|\De G\|_{\infty}$. 
      
To study the global attractor, we need to derive 
some a priori estimates about solutions. 
 
\begin{lemma} 
 Suppose $u$ is a solution of (\ref{sh4})-(\ref{sh6}).
 Then  $u$ is uniformly (in time) 
 bounded,  and  the following estimates hold for $t>0$
 
 \be 
 \|u(x,y,t)\|^2 \leq \|u_0(x,y)\|^2 exp(-2\mu t)+ \frac{\mu}{b}, \; 
  \label{sh8} 
 \ee 
 and thus
\be 
 {\lim \sup}_{t\rightarrow +\infty} \|u(x,y,t)\|    \leq  \sqrt{\frac{\mu}{b}} \equiv R ,   
\label{sh9} 
\ee 
where $R = \sqrt{\frac{\mu}{b}}$.
\label{bound} 
 
\end{lemma}

{\bf Proof.}  Taking the inner product of (\ref{sh4}) with $u$, we have 
 
\be 
\frac{1}{2}\frac{d}{dt}\|u\|^2 + \|\De u\|^2 + 2(\De u,u) + \a \|u\|^2 
				\nonumber  \\
  + (u^2,\id \sqg u^2 (\xi,\eta)d\xi d\eta)   = 0. 
\label{sh10} 
\ee 
Note that 
$$  2|(\De u,u)|\leq 2\|\De u\|\|u\|  \leq  \|\De u\|^2 + \|u\|^2,$$ 
 
$$  (u^2,\id \sqg u^2(\xi,\eta)d\xi d\eta) $$ 
$$ = \id u^2(\id \sqg u^2(\xi,\eta)d\xi d\eta)dxdy$$ 
 
$$ \geq b\id u^2(x,y)dxdy \id u^2(\xi,\eta)d\xi d\eta = b\|u\|^4. 
$$ 
Then from (\ref{sh10}) we get 
\be
 \frac{d}{dt}\|u\|^2 + 2(\a-1)\|u\|^2 + 2b\|u\|^4 \leq 0. 
 \label{L2norm}
 \ee
It is easy to see that if $\a \geq 1$, i.e., $\mu \leq 0$,
then  all solutions
approach zero in $L^2$. We will not consider this simple dynamical
case. In the rest of this paper we assume that $\mu > 0$, i.e.,
$\a < 1$.
 
Thus we have, for any constant $\e >0$,
\be
\frac{d}{dt}\|u\|^2 +  2\e \|u\|^2 +2(\a-1-\e)\|u\|^2+ 2b\|u\|^4 \leq  0,
\ee 
or
\be
\frac{d}{dt}\|u\|^2 +  2\e \|u\|^2 
 + [ \frac{(\a-1-\e)}{\sqrt{2b}} + \sqrt{2b}\|u\|^2 ]^2  
 	\leq \frac{(\a-1-\e)^2}{2b}.
\ee  
 So
\be
\frac{d}{dt}\|u\|^2 +  2\e \|u\|^2  \leq  \frac{(\a-1-\e)^2}{2b}.
\ee  
By the usual Gronwall inequality (\cite{Temam}) we obtain
\be
 \|u\|^2 \leq \|u_0\|^2 exp(-2\e t)+ \frac{(\a-1-\e)^2}{4b\e}.
\ee
When $\e = 1-\a = \mu$, we get the optimal or tight estimate
\be
 \|u\|^2 \leq \|u_0\|^2 exp(-2\mu t) + \frac{\mu}{b}.
\ee
This completes the proof of  
Lemma \ref{bound}.   $\hfill \qed$

Moreover, higher order derivatives of $u$ are also
uniformly bounded. 
 
\begin{lemma} 
 Suppose $u$ is a solution of (\ref{sh4})-(\ref{sh6}). Then $\|\na u\|$ 
 and $\|\De u\|$ are uniformly (in time)  bounded. 
\label{highbound} 
\end{lemma}

In order to prove this lemma, we recall a few useful inequalities.

{\em Uniform Gronwall inequality} (\cite{Temam}).  
Let $g,h,y$ be three positive locally  
integrable functions on $[t_0,+\infty)$ satisfying the inequalities 
$$  
\frac{dy}{dt}\leq gy + h , 
$$ 
with  $\int_t^{t+1}gds \leq a_1,$  $\int_t^{t+1}hds\leq a_2$ 
and $\int_t^{t+1}yds \leq a_3$ for $t\geq t_0,$ 
where the $a_i$(i=1,2,3) are positive constants. Then 
$$  
y(t+1) \leq (a_2+a_3)exp(a_1), \mbox{for} \; t \geq t_0 . 
$$

{\em Gagliardo-Nirenberg inequality} (\cite{Pazy}). 
Let $w \in L^q\cap W^{m,r}(D)$, where $1 \leq q,r \leq \infty$.  
For any integer $j$, $0\leq j \leq m$, $\frac{j}{m}\leq \la \leq 1.$ 
$$  
\|D^j w\|_p \leq C_0 \|w\|_q^{1-\la}\|D^m w\|_r^{\la}  
\label{GN} 
$$  
provided 
$$\frac{1}{p} = \frac{j}{n} + \la (\frac{1}{r} - \frac{m}{n})  
+ \frac{1-\la}{q},$$ 
and $m-j-\frac{n}{r}$ is not a nonnegative integer 
If $m-j-\frac{n}{r}$ is a nonnegative integer, then the inequality  
(\ref{GN}) holds for  
$\la = \frac{j}{m}$.   
 
{\em Poincar\'e inequality} (\cite{Friedman}). 
For $w \in H_0^1(D)$, 
$$ \la_1 \| w \|^2 \leq \|\na w \|^2 , $$ 
where $\la_1$ is the first eigenvalue of $-\De$ on the domain $D$, with zero  
Dirichlet boundary condition  
on $\p D$.

{\bf Proof of Lemma \ref{highbound}}.   
Due to the boundary condition (\ref{sh6}) on $ \na u$ 
and the Poincar\'e  inequality, we get  
$ \| \na u \|^2 \leq \la_1^{-1} \|\De u \|^2$. Hence it is  
sufficient to prove that $ \|\De u\|$ is 
bounded.   We first show that $\int_t^{t+1}\|\De u\|^2ds$ 
is bounded. In fact, using 
 $$  2|(\De u,u)| \leq 2\|\De u\|\|u\|  \leq \frac12 \|\De u\|^2 + 2\|u\|^2,$$
in (\ref{sh10}), we get
\be
 \frac{d}{dt} \|u\|^2 + \|\De u\|^2 + 2(\a-2)\|u\|^2 + 2b\|u\|^4 \leq 0. 
\ee
Since 
$$2b\|u\|^4 + 2(\a-2)\|u\|^2 = b\|u\|^2 +  
2\b(\|u\|^4 + \frac{2\a-4-\b}{2\b}\|u\|^2)$$
$$=b\|u\|^2 + 2\b(
\|u\|^2 + \frac{2\a-4-\b}{4\b})^2 - 
\frac{(2\a -4 -\b)^2}{8\b}$$
$$\geq b\|u\|^2 - \frac{(2\a -4 -\b)^2}{8\b},$$
we conclude
\be
\frac{d}{dt} \|u\|^2 + \|\De u\|^2 + b\|u\|^2 \leq 
\frac{(2\a -4 -\b)^2}{8\b} =\frac{(2 + 2\mu + \b)^2}{8\b}.
\label{sh11.5}
\ee
Integrating (\ref{sh11.5}) with respect to $t$ from $t$ to $t+1$ and noting
Lemma \ref{bound}, we see that $\int_t^{t+1}\|\De u\|^2ds$ is bounded.

Now, multiplying (\ref{sh4}) by $\De^2 u$ and integrating over $D$, 
it follows that 
 
$$ 
\frac{1}{2}\frac{d}{dt}\|\De u\|^2 + \|\De^2 u\|^2 +  
2\id \De u\De^2udxdy + \a \|\De u\|^2$$ 
\be 
+\id (u \id \sqg u^2(\xi,\eta)d\xi d\eta)\De^2 udxdy = 0. 
\label{sh12} 
\ee 
Note that 
\be 
2|\id\De u\De^2 udxdy |  \leq \frac{1}{2}\|\De^2u\|^2 + 2\|\De u\|^2,  
\ee 
and 
$$ 
 | \id ( u \id \sqg u^2(\xi,\eta)d\xi d\eta)  \De^2  udxdy |  
$$ 
$$ 
  =   | \id (\De u)^2(\id \sqg  u^2(\xi,\eta)d\xi d\eta)dxdy  
$$  
$$ 
       +   \id u\De u (\id(\De \sqg))u^2(\xi,\eta)d\xi d\eta)dxdy  
$$ 
$$  
   +   2\id \na u\De u(\id(\na \sqg ))u^2(\xi,\eta)d\xi d\eta)dxdy|   
$$  
$$ 
   \leq    (a\|u\|^2+2\la_1^{-\frac12} \|\na G \|_{\infty}\|u\|^2
   + \frac12 \|\De G \|_{\infty}\|u\|^2 ) \|\De u\|^2 
	       + \frac12 \|\De G \|_{\infty} \|u\|^4 
$$  
$$          
   \leq    (a +2\la_1^{-\frac12} K_1  + \frac12 K_2 )\|u\|^2 \| \De u \|^2 
   			+ \frac12 K_2  \|u\|^4,   
$$ 
where $a, K_1, K_2$ are various upper bounds of $G$ defined
in (\ref{sh7}), and $R$ is the $L^2$ bound of the
solution $u$ as in Lemma 1.
Hence by (\ref{sh12}) we get  
\be 
\frac{d}{dt}\|\De u\|^2 \leq 2 [(a +2\la_1^{-\frac12} K_1 + \frac12 K_2)\|u\|^2
	 -\a +2]  \|\De u\|^2   +  K_2  \|u\|^4. 
\label{sh13} 
\ee 
Finally, applying the uniform Gronwall inequality  (\ref{sh13}) 
and noting Lemma \ref{bound}, we conclude 
 that $\| \De u\|^2$ is uniformly  
bounded for all $t\geq 0.$ 
This proves Lemma \ref{highbound}.  $\hfill \qed$

We now have the following global existence and uniqueness 
result. 
 
\begin{theorem} 
Let $u_0(x,y) \in L^2(D)$ and $G$ satisfies (\ref{sh7}),  
then the initial-boundary value problem 
$(\ref{sh}), (\ref{sh5}), (\ref{sh6})$ has a unique global 
solution $u\in L^{\infty}(0,\infty; H^2_0(D))$.  Moreover, 
the corresponding solution semigroup $S(t)$, defined by 
$$ 
u = S(t) u_0, 
$$ 
has a bounded absorbing set   
$$ 
B_0 = \{ u  \in H^2_0(D):  (\|u\|^2+\|\na u\|^2 
   +\|\De u\|^2)^{\frac{1}{2}}  \leq  \tilde{R}  \},
$$ 
where $\tilde{R} $ is a postive constant which depending on the  
uniform  bound of $\|u\|, \|\na u\|, \|\De u\|$. 
Finally, the solution semigroup $S(t)$, when restricted on 
$H^2_0(D)$, is continuous from  
$H^2_0(D)$ into $H^2_0(D)$ for $t > 0$. 
\label{global} 
\end{theorem}

{\bf Proof.}  
The global existence, uniqueness and absorbing property
follow from standard arguments (e.g., \cite{Temam})
together with Lemmas \ref{bound}, \ref{highbound} above.  
The absorbing property also follows from these two lemmas. 
						 
We now prove that $S(t)$ is continuous in $H^2(D)\cap H_0^1(D)$.                                              
Suppose that $u_0,v_0\in H^2_0(D)$ with $\|\De u_0\|, 
\|\De v_0\|\leq 2R_1,$ 
we denote by $u(t), v(t)$ the corresponding solutions, i.e., 
$u(t) = S(t)u_0, v(t) = S(t)v_0$.  Let  
$w(t)=u(t)-v(t).$ Then $w(t)$ satisfies 
$$ 
w_t + \De^2 w + 2 \De w + \a w + w\id \sqg u^2(\xi,\eta) 
d\xi d\eta+$$ 
\be 
v\id \sqg (u(\xi,\eta)+v(\xi,\eta))w(\xi,\eta)d\xi d\eta = 0. 
\label{sh14} 
\ee 
Applying the Gagliardo-Nirenberg inequality 
 
 $$ \|u\|_{\infty} \leq C_0\|\De u\|,$$ 
and the Poincar\'e inequality   
 
$$ \|w\| \leq \frac{1}{\la_1}\|\De w\|,$$ 
we obtain (similar to the proof of Lemma \ref{highbound}),   
 
$$  
\frac{d}{dt}\|\De w\|^2 \leq C_1 \|\De w\|^2, 
$$  
which implies that $\|\De w(t)\|^2 \leq \|\De w_0\|^2 exp(C_1 t)$ 
for some positive constant $C_1$. This shows that $S(t)$ is continuous.                                               
						$\hfill \qed$

This theorem implies that  (\ref{sh4})-(\ref{sh6})
defines an infinite dimensional nonlocal dynamical system.

In the rest of this section, we consider the global attractor
for the nonlocal dynamical system (\ref{sh4})-(\ref{sh6}). 
We will establish the following result about the global attractor. 
 
\begin{theorem}
There exists a global attractor ${\cal A}$
for the nonlocal dynamical system 
$(\ref{sh}), (\ref{sh5}), (\ref{sh6})$.
The global attractor is the $\omega-$limit set 
of the absorbing set $B_0$ (as in Theorem \ref{global}), and it   
has the following properties:
 
(i) $A$ is compact and $S(t){\cal A} = {\cal A},$ for $t > 0$;  
 
(ii) for every bounded set $B\subset H^2_0(D)$,  
$\lim\limits_{t\to \infty} d(S(t)B,{\cal A})=0;$ 
 
(iii)${\cal A}$ is connected in $H^2_0(D),$  
where  
$d(X,Y)=\sup\limits_{x\in X}\inf\limits_{y\in  Y}\|x-y\|_{H^2_0(D)}$  
is the Hausdorff distance. 

Moreover, the global attractor ${\cal A}$ has finite Hausdorff dimension  
$d_H({\cal A}) \leq m$,  
where   
$$ 
m \sim  C (1+\sqrt{ \mu + (2a-b)\frac{\mu}{b} }), 
\label{attractor} 
$$ 
where $ C >0$ is a constant depending only on the domain $D$,
and $a>0, b>0$ are the upper, lower bounds of the kernel $G$,
respectively.
\end{theorem}
 
{\bf Proof.}  
The existence and properties of ${\cal A}$ are 
quite standard now (see \cite{Temam} and references therein). 
We omit this part, and only estimate the dimensions below.

As in  \cite{Temam}, we may use 
the so-called Constantin-Foias-Temam trace formula (which works for 
the semiflow $S(t)$ here) to estimate 
the sum of the global Lyapunov exponents of 
${\cal A}$. The sum of these Lyapunov exponents can then 
 be used to estimate the upper bounds of ${\cal A}$'s 
Hausdorff dimension, $d_H({\cal A})$. 
To this end,  
we linearize equation (\ref{sh4}) about a solution $u(t)$ in the  
global attractor to obtain an equation for $v(t)$ and then
use the trace formula to estimate the sum of 
the global Lyapunov exponents. Doing so, we obtain
\be 
v_t + L(u(t))  v = 0, 
\label{sh15} 
\ee 
where    
$$ 
L(u(t))v = \De^2 v + 2\De v + \a v +  
v\id \sqg u^2(\xi,\eta)d\xi d\eta  
$$ 
$$+2u\id \sqg u(\xi,\eta) v(\xi,\eta)d\xi d\eta. 
$$ 
This equation is supplemented with 
$v(x,y , 0) = \xi(x,y) \in  H^2_0(D)$.
Denote by $\xi_1(x,y), \ldots , \xi_m(x,y)$, 
$m$ linearly independent functions
in $H^2_0(D)$, and $v_{i}(x,y , t)$ the solution of (\ref{sh15})
satisfying $v_{i}(x,y , 0) = \xi_{i}(x,y)$, $i = 1, \ldots, m$.
Let $Q_m (t)$ represent the orthogonal projection of $H^2_0(D) $
onto the subspace spanned by $\{v_1(x,y,t), \ldots , v_m(x,y, t) \}$. 

We need to estimate the lower bound of 
$Tr(L(u(t)Q_m(t)))$, which gives bounds on the sum of 
global Lyapunov exponents. 
Note that in \cite{Temam},
the linearized equation like (\ref{sh15}) is written
as $v_t = L(u(t))v$ and in that case one needs to
estimate the upper bound of $Tr(L(u(t)Q_m(t)))$.
Suppose that $\phi_1(t),...,\phi_m(t)$ is   
an orthonormal basis  ($\|\phi_j\| =1$)
of the subspace $Q_m(t) H_0^2(D)$
for any $t>0$. 

Now we estimate  the lower bound of 
$Tr(L(u(t)Q_m(t)))$. It is easy to see that
 
$$Tr(L(u(t)Q_m))$$ 
 
$$ 
=\sum\limits^m_{j=1}(\De^2\phi_j+2\De\phi_j+\a \phi_j+\phi_j\id 
	\sqg u^2(\xi,\eta)d\xi d\eta,\phi_j) 
$$ 
$$ 
+ \sum\limits^m_{j=1}(2u\id \sqg u(\xi,\eta) \phi_j d\xi d\eta,\phi_j). 
$$ 
Since 
$(2\De\phi_j, \phi_j) \geq -(\frac1{\e}\|\De \phi_j\|^2 +\e \|\phi_j\|^2)$ 
for any constant $\e > 1$, we get 
$$ 
Tr(L(u(t)Q_m)) \geq \sum\limits^m_{j=1} [(1-\frac1{\e})\|\De \phi_j\|^2+  
b\|\phi_j\|^2\|u\|^2 + (\a - \e) \|\phi_j\|^2 ] 
$$ 
$$ 
+ \sum\limits^m_{j=1} 2\id u\phi_jdxdy (\id \sqg u(\xi,\eta) \phi_j d\xi d\eta) 
$$ 
$$ 
\geq \sum\limits^m_{j=1} (1-\frac1{\e}) \|\De \phi_j\|^2 
   + \sum\limits^m_{j=1} (b\|u\|^2+ \a-\e -2 a\|u\|^2)        
$$ 
\be 
 = \sum\limits^m_{j=1} (1-\frac1{\e}) \|\De \phi_j\|^2 
  + [1-\mu-\e + (b-2 a) \|u\|^2 ] m. 
\label{sh18} 
\ee 

We introduce notation $f(x,y) = \sum\limits_{j=1}^m | \phi_j|^2$.
Note that $m = \id f(x,y)dxdy$.
By the generalized Sobolev-Lieb-Thirring  
inequality (\cite{Temam}, page 462), 
 
$$\id f^3(x,y)dxdy \leq K_0 \sum\limits^m_{j=1} \|\De \phi_j\|^2,$$ 
where $K_0>0$ depending only on the domain $D$. 
Moreover, due to the fact that $L^3 (D) \hookrightarrow L^1(D)$,
$$m^3 = (\id f(x,y)dxdy)^3 \leq C_2 \id f^3(x,y)dxdy $$ 
 
$$\leq K_0 C_2 \sum\limits^m_{j=1} \|\De \phi_j\|^2 $$ 

$$ = C\sum\limits^m_{j=1} \|\De \phi_j\|^2 $$ 
for some constants $C_2>0, C>0$ depending only on the domain $D$.
 
Thus 
\be 
(1-\frac1{\e}) \sum\limits^m_{j=1} \|\De \phi_j\|^2
  \geq    (1-\frac1{\e}) \frac1{C} m^3.
     \label{sh19} 
\ee

Therefore, by (\ref{sh18})-(\ref{sh19}) we have 
\be 
Tr (L(u(t)Q_m))
& \geq & \frac{1-\frac1{\e}}{C} m^3 - (\mu-1 +\e + (2a-b) \|u\|^2) m
		 		\nonumber   \\
& \geq & \frac{1-\frac1{\e}}{C} m^3 - (\mu-1 +\e + (2a-b)\frac{\mu}{b}) m
				\nonumber   \\
& > & 0
\ee
whenever 
\be
m >\sqrt{[ \mu-1 +\e + (2a-b) \frac{\mu}{b} ]
	\frac{C}{1-\frac1{\e}}}.
\label{mm}
\ee
The right hand side of (\ref{mm}) has the minimal value of 
\be
   m \sim    C (1+\sqrt{\mu + (2a-b)\frac{\mu}{b}})
   \label{dimension}
\ee
when 
$ \e = 1+\sqrt{ \mu + (2a-b)\frac{\mu}{b} }$.
 
As in \cite{Temam},  we conclude that the Hausdorff dimension of  
${\cal A}$ is  estimated as in  (\ref{dimension}).
This proves Theorem \ref{attractor}.       $\hfill \qed$

 \section{Local Swift-Hohenberg Model}

Similarly, for the two-dimensional local Swift-Hohenberg 
equation (\ref{sh3}), we can obtain the existence of 
the global attractor $\tilde {\cal A}$. We omit this part
and will only estimate the dimension 
of $\tilde {\cal A}$.  

\begin{theorem} 
There exists the  global attractor $\tilde {\cal A}$  for the 
local dynamical system (\ref{sh3}), (\ref{sh5}), (\ref{sh6}). 
The Hausdorff dimension  of  $\tilde {\cal A}$ is finite, and 
$d_H(\tilde {\cal A} )  \leq  m_1 \sim  C (1+\sqrt{ \mu }) $, 
where $ C $ is a constant depending only on the domain $D$. 
 \label{attractor2}  
\end{theorem} 
  
{\bf Proof.} As in the proof of Ttheorem \ref{attractor},
 we consider the linearized equation   of (\ref{sh3}), defined by 
$$ v_t + L_1(u(t))v = 0,$$
where
$$L_1(u(t))v = \De^2v + 2\De v + \a v + 3u^2v.  $$ 

Then we estimate
$$ Tr(L_1(u(t)Q_m)) $$ 
 
$$ 
=\sum\limits^m_{j=1}(\De^2\phi_j+2\De\phi_j+\a \phi_j+3u^2\phi_j ,\phi_j) 
$$ 
$$
= \sum\limits^m_{j=1}[\|\De \phi_j\|^2 + 2(\De \phi_j, \phi_j)
	+\a \|\phi_j\|^2 + 3(u^2\phi_j, \phi_j)]
$$
$$ 
\geq \sum\limits^m_{j=1} (1-\frac1{\e})\|\De \phi_j\|^2  
  + \sum\limits^m_{j=1}(\a -\e), 
$$ 
where we have used the fact that  $3(u^2\phi_j, \phi_j) \geq 0$.
Noting   again that
$m^3 \leq  C\sum\limits^m_{j=1} \|\De \phi_j\|^2$ and $\a =1-\mu$,
we have
\be
Tr (L_1(u(t)Q_m))
\geq \frac{1-\frac1{\e}}{C} m^3 -(\mu -1+\e) m >0
\ee
whenever 
\be
m >\sqrt{( \mu-1 +\e ) \frac{C}{1-\frac1{\e}} }.
\label{mmm}
\ee
The right hand side of (\ref{mmm}) has the minimal value of 
\be
   m \sim  C (1+\sqrt{ \mu })
\ee
when 
$ \e = 1+\sqrt{ \mu }$.  
This completes the proof.      		$\hfill \qed$

\section{Discussions}

In this paper, we have discussed the Hausdorff dimension estimates 
for the global attractors of the two-dimensional nonlocal and local  
Swift-Hohenberg model for Rayleigh-Benard convection. 

The Hausdorff dimension for the global attractor of the
nonlocal model is estimated as
$$
   m \sim    C (1+\sqrt{\mu + (2a-b)\frac{\mu}{b}}),
$$
while for the local model this estimate is
$$
   m \sim  C (1+\sqrt{ \mu }),
$$
where $C>0$ is an absolute constant depending only on
the fluid convection domain, and $\mu>0$ measures the 
difference of the Rayleigh number from its critical 
convection onset value.
Note that $a, b > 0$ are the upper and lower
bounds, respectively, of the kernel $G$ of the nonlocal nonlinearity
in (\ref{sh}).

The two dimension estimates above
differ by an absolute constant $(2a-b)\frac{\mu}{b}$, 
which depends only on
the the Rayleigh number through $\mu$, and upper and 
lower bounds of the kernel $G$ of the nonlocal
nonlinearity.  Moreover, if the kernel $G$ is a constant function
(thus, $a=b=G$), then the dimension estimate for the
nonlocal model becomes 
$$
   m \sim    C (1+\sqrt{ 2 \mu}),
$$
which still differs from the dimension
estimate for the local model by a constant depending on the Rayleigh number
through $\mu$.

\bigskip
 
{\bf Acknowledgement.} 
Part of this work was done while Jinqiao Duan was visiting the
Institute for Mathematics and its Applications (IMA), Minnesota,
and the Center for Nonlinear Studies, Los Alamos National
Laboratory.
    This work was supported by the Nonlinear Science Program of China,   
the National Natural Science Foundation of China Grant 19701023, 
the Science Foundation of Chinese  
Academy of Engineering Physics  Grant 970682, and the USA National Science 
Foundation Grant DMS-9704345.

\end{document}